\documentclass[11pt,a4paper]{article}

\usepackage{amsmath, amssymb} 

\newtheorem{corollary}{Corollary}
\newtheorem{definition}{Definition}
\newtheorem{example}{Example}
\newtheorem{remark}{Remark}
\newtheorem{lemma}{Lemma}
\newtheorem{proposition}{Proposition}
\newtheorem{theorem}{Theorem}

\def\<{\langle}
\def\>{\rangle}
\def\RR{\mathbb{R}}
\newcommand\tr{\operatorname{Tr}}
\def\Ric{\operatorname{Ric}}

\author{Vladimir Rovenski\footnote{Mathematical Department, University of Haifa, Mount Carmel, 3498838 Haifa,  Israel
        \newline e-mail: {\tt vrovenski@univ.haifa.ac.il}        }
       }

\title{Weak Almost Contact Structures: a Survey}

\begin{document}

\date{}

\maketitle

\begin{abstract}
Weak almost contact structures,
i.e., the~linear complex structure on the contact distribu\-tion is replaced by a nonsingular skew-symmetric tensor,
defined by the author and R.\,Wolak, allowed us to take a new look at the theory of contact manifolds.
The paper surveys recent results (concerning geodesic and Killing fields, rigidity and splitting theorems,
Ricci-type solitons and Einstein-type metrics, warped products, etc.) in this new field of Riemannian~geometry.

\vskip1mm\noindent
\textbf{Keywords}:
almost contact structure, K-contact structure,
(nearly) Sasakian and cosymplectic structures, $\beta$-Kenmotsu manifold, Ricci-type soliton,
Einstein-type metric

\textbf{Mathematics Subject Classifications (2010)} Primary 53C12; Secondary 53C21
\end{abstract}

\section{Introduction}
\label{sec:01}
 \vglue-10pt
 \indent

{Contact Riemannian geomet\-ry} is of growing interest due to its important role in mechanics, \cite{blair2010riemannian,OV-2024}.
Sasakian manifolds (normal contact metric manifolds) and, more generally, $K$-contact manifolds
have become an important and active subject, especially after the appearance of the fundamental treatise~\cite{bg2} of C.\,Boyer and K.\,Galicki.
On a ${K}$-contact manifold $M(f,\xi,\eta,g)$ the vector field $\xi$ is Killing (the Lie derivative $\pounds_\xi\,g=0$),
and geodesic ($\nabla_{\xi}\,\xi=0$).
The~restriction of $f$ to $f(TM)$ determines a complex structure.
If~a plane in $TM$ contains $\xi$, then its sectional curvature is called $\xi$-\textit{sectional curvature}.
The $\xi$-sectional curvature of a $K$-contact manifold~is equal to one.
Two~important subclasses of ${K}$-contact manifolds are
\begin{equation}\label{E-nS-Sas}
 (\nabla_X f)Y =
 \left\{\begin{array}{cc}
  g(X,Y)\,\xi -\eta(Y)X\,, & {\rm Sasakian,} \\
  0\,, & {\rm cosymplectic.}
 \end{array}\right.
\end{equation}
Any cosymplectic manifold is locally the product of a K\"{a}hler manifold and $\RR$.
A~Riemannian manifold $(M^{2n+1},g)$ with a {contact} 1-form $\eta$ (i.e., $\eta\wedge(d\eta)^n\ne0$)
is Sasakian, if its
cone $M\times\mathbb{R}^{>0}$ with the
metric $t^2 g+dt^2 $ is a K\"{a}hler manifold.
We~get $\nabla_{X}\,\xi=-\frac12\,fX$ on Sasakian manifolds, and $\nabla_{X}\,\xi=0$ on cosymplectic manifolds.

Nearly K\"{a}hler {mani\-folds} \cite{G-70} are defined by the condition that the symmetric part of $\nabla J$ vanishes,
in contrast to the K\"{a}hler case where $\nabla J=0$.
Nearly Sasakian/cosymp\-lectic manifolds are defined similarly
-- by a constraint on the symmetric part of the structural tensor $f$
-- starting from Sasakian/cosymplectic manifolds:
\begin{equation}\label{E-nS-01}
 (\nabla_X{f})Y + (\nabla_Y{f})X =
 \left\{\begin{array}{cc}
  \hskip-1mm 2\,g(X,Y)\xi -\eta(Y)X -\eta(X)Y\,, & {\rm nearly\ Sasakian,} \\
  0\,, & \hskip-6mm {\rm nearly\ cosymplectic.}
 \end{array}\right.
\end{equation}
In dimensions greater than 5, every nearly Sasakian manifold is Sasakian,
and a nearly cosymplectic manifold $M^{2n+1}$ splits into ${\mathbb R}\times F^{2n}$ or $B^5\times F^{2n-4}$,
where $F$ is a nearly K\"{a}hler manifold and $B$ is a 5-dimensional nearly cosymplectic manifold.
These structures
appeared in the study of harmonic almost contact manifolds.
For nearly Sasakian/cosymplectic manifolds we get $g(R_{\,\xi, Z}\,{f} X, {f} Y) = 0$, see \cite{E-2005,NDY-2018},
where $R_{{X},{Y}}Z=\nabla_X\nabla_Y Z -\nabla_Y\nabla_X Z -\nabla_{[X,Y]} Z$ is the Riemann tensor;
thus,
\begin{equation}\label{E-nS-04cc}
 R_{X,Y}Z\perp\xi\quad (X,Y,Z\perp\xi) -
 \mbox{the distribution $\ker\eta$ is \textit{curvature invariant}}.
\end{equation}
For example,
the distribution $\ker\eta$ of any 1-form $\eta$ on $\RR^m$
 satisfies~\eqref{E-nS-04cc}.

The~Ricci tensor is given by ${\rm Ric}({X},{Y})=\tr(Z\to R_{Z,X}Y)=\!\sum_{i} g(R_{e_i,{X}}{Y},e_i)$,
where $(e_i)$ is any local orthonormal basis of~$TM$.
The Ricci operator $\Ric^\sharp$
associated with the Ricci tensor
is given by
 ${\rm Ric}({X},{Y})=g(\Ric^\sharp{X},{Y})$.
The scalar curvature of $(M,g)$ is given by $r={\rm trace}_g \Ric$.
On some compact manifolds there are no Einstein metrics, which motivates the study of generalizations of such metrics.
The~\textit{genera\-lized Ricci soliton}
is given for some smooth vector field ${V}$ and real $c_1, c_2$ and $\lambda$~by
\begin{equation}\label{E-g-r-e}
 (1/2)\,\pounds_{{V}}\,g = -c_1 {V}^\flat\otimes {V}^\flat + c_2\Ric + \lambda\,g.
\end{equation}
If ${V} = \nabla f$ in \eqref{E-g-r-e} for some $f\in C^\infty(M)$, then by
the definition
${\rm Hess}_f({X},{Y})=\frac12\,(\pounds_{\nabla f}\,g)({X},{Y})$,
we get the \textit{generalized gradient Ricci soliton} equation, see  \cite{CZB-2022}:
\begin{equation}\label{E-gg-r-e2}
 {\rm Hess}_{f_1} = -c_1 df_2\otimes df_2 + c_2\Ric + \lambda\,g
\end{equation}
for some
$f_1,f_2\in C^\infty(M)$ and real $c_1, c_2$ and $\lambda$.
For different values of $c_1,\,c_2,\lambda$, equation \eqref{E-g-r-e} is a generalization of
Einstein metric, $\Ric + \lambda\,g=0$ ($c_1=0,\,c_2=-1,\, V=0$),
Killing equation ($c_1 = c_2 = \lambda = 0$),
Ricci soliton equation
($c_1 = 0,\,c_2 = -1$), etc.,
see~\cite{G-D-2020}.
In~\cite{cho2009ricci}, Cho-Kimura
 defined $\eta$-\textit{Ricci soliton} on a
contact metric manifold $M^{2n+1}({f},\xi,\eta,g)$ by the following equation:
\begin{equation}\label{1.1}
 ({1}/{2})\,\pounds_{V}\,g +{\rm Ric} +\lambda\,g +\mu\,\eta \otimes\eta =0,
\end{equation}
where
$\lambda,\mu\in\RR$.
For $\mu=0$, \eqref{1.1} gives a Ricci soliton; moreover, if $V=0$, then \eqref{1.1} gives an Einstein manifold.
 A contact metric manifold
 is called $\eta$-\textit{Einstein}, if
 ${\rm Ric}= a\,g + b\,\eta\otimes\eta$
is true,
where $a,b\in C^\infty(M)$,
see~\cite[p.~285]{YK-1985} for $a,b\in\mathbb{R}$.
Note that $a$ and $b$
can be non-constant for
an $\eta$-Einstein Kenmotsu manifold,
see~\cite{kenmotsu1972class}.

Many articles study the questions: How interesting Ricci-type solitons are for contact metric mani\-folds?
When an almost contact metric manifold
equipped with a Ricci-type soliton
carries Einstein-type~metrics?

In \cite{RWo-2,rov-119}, we introduced metric structures on a smooth manifold that generalize the almost contact, Sasakian, cosymplectic, etc. metric structures.
 Such so-called ``weak" structures
(the complex structure on the contact distribution $f(TM)$ is replaced by a nonsingular skew-symmetric tensor)
allow us to take a fresh look at the theory of classical structures and find new applications.
A \textit{weak almost contact structure} on a smooth manifold $M^{2n+1}$ is defined by a $(1,1)$-tensor field $f$ of rank $2n$,
a vector field $\xi$, a 1-form $\eta$ and a~nonsingular $(1,1)$-tensor field $Q$ satisfying, see \cite{RWo-2},
\begin{equation}\label{E-nS-2.1}
 f^2 = -Q +\eta\otimes\xi, \qquad \eta(\xi)=1,\qquad Q\,\xi=\xi .
\end{equation}
A ``small" (1,1)-tensor $\widetilde Q=Q-{\rm id}$ measures the difference between a weak almost contact structure and an almost contact one.
We assume that a $2n$-dimensional distribution $\ker\eta$ is ${f}$-invariant (as in the classical theory \cite{blair2010riemannian}, where $Q={\rm id}$).

If~there exists a Riemannian metric $g$ such that
\begin{equation}\label{2.2}
 g({f} X,{f} Y)= g(X,Q\,Y) - \eta(X)\,\eta(Y),\quad X,Y\in\mathfrak{X}_M,
\end{equation}
then $({f},Q,\xi,\eta,g)$ is called a {\it weak almost contact metric structure},
and $g$ is a \textit{compatible metric}.
 Putting $Y=\xi$ in \eqref{2.2} and using $Q\,\xi=\xi$, we get
 $g(X,\xi) = \eta(X)$;
moreover,
the tensor ${f}$ is skew-symmetric and the tensor $Q$ is self-adjoint.

\begin{remark}\rm
(i)~The concept of an almost para-contact structure
is closely related to an almost contact structure and an almost product structure, see \cite{FP-2017}.
Similarly to \eqref{E-nS-2.1}, we define a \textit{weak almost para-contact structure} by
$f^2 = Q -\eta\otimes\xi,\ Q\,\xi=\xi$, see \cite{RWo-2}.

(ii)~A weak almost contact structure admits a compatible Riemannian metric if ${f}$
admits a skew-symmetric representation, i.e., for any $x\in M$ there exist a neighborhood $U_x\subset M$ and a~frame $\{e_k\}$ on $U_x$,
for which ${f}$ has a skew-symmetric matrix.
\end{remark}

A~\textit{weak contact metric structure} is a weak almost contact metric structure satis\-fying $d\eta = \Phi$,
where the {$2$-form} $\Phi$ is defined by $\Phi({X},{Y})=g({X},{f} {Y}),\ {X},{Y}\in\mathfrak{X}_M$.
For~a weak contact metric structure $({f},Q,\xi,\eta,g)$, the 1-form $\eta$ is contact.
A weak almost contact structure is \textit{normal} if the tensor ${\cal N}^{\,(1)}= [{f},{f}] + 2\,d\eta\otimes\,\xi$ vanishes.

A weak almost contact metric structure
is called a \textit{weak ${K}$-struc\-ture} if it is normal and
$d \Phi=0$.
 We define two subclasses of weak ${K} $-manifolds as follows:
\textit{weak cosymplectic manifolds} if $d\eta=0$, and \textit{weak Sasakian manifolds}~if
 $d\eta = \Phi$.
Omitting the normality condition, we get the following:
a weak almost contact metric structure
is called
a \textit{weak almost cosymplectic structure} if $d\eta = \Phi$ is valid;
and
a \textit{weak almost Sasakian structure} if $\Phi$ and $\eta$ are closed forms.

In~\cite{rov-119}, we proved rigidity results: a weak Sasakian structure is the Sasakian structure
and a weak almost contact metric structure satisfying $\nabla f=0$ is a weak cosymplectic structure.
For this, we calculated the derivative $\nabla f$, using a new tensor ${\cal N}^{(5)}$, see \eqref{E-N5},
in addition to classical tensors ${\cal N}^{(i)}\ (i=1,2,3,4)$.

\begin{example}\rm
Let $M^{2n+1}(\varphi,Q,\xi,\eta)$ be a weak almost contact manifold.
Consider the product manifold $\bar M = M^{2n+1}\times\mathbb{R}$,
and define tensor fields $\bar f$ and $\bar Q$ on $\bar M$ putting
\begin{equation*}
 \bar f(X,\, a\partial_1) = (fX- a\,\xi,\, \eta(X)\partial_1),\quad
 \bar Q(X,\,  a\,\partial_1) = (QX,\, a\partial_1),\quad a\in C^\infty(M).
\end{equation*}
Hence, $\bar f(X,0)=(fX,0)$, $\bar Q(X,0)=(QX,0)$ for $X\in\ker f$,
$\bar f(\xi,0)=(0,\partial_1)$, $\bar Q(\xi,0)=(\xi,0)$
and
$\bar f(0,\partial_1)=(-\xi,0)$, $\bar Q(0,\partial_1)=(0,\partial_1)$.
Then it is easy to verify that $\bar f^{\,2}=-\bar Q$. The tensors ${\cal N}^{\,(i)}\ (i=1,2,3,4)$ appear when we use
the integrability condition $[\bar f, \bar f]=0$ of $\bar f$ to express the normality condition ${\cal N}^{\,(1)}=0$ of a weak almost contact structure.
\end{example}

The~relationships between some classes of weak structures is summarized in the~following diagram (known for classical~structures):
\[
\left|
   \begin{array}{c}\hskip-2mm
  \textrm{weak\ almost} \\
  \textrm{contact\ metric}
   \hskip-1mm\end{array}
 \right|
\overset{\Phi=d\eta}\longrightarrow
  \left|
   \begin{array}{c}\hskip-1mm
  \textrm{weak\ contact} \\
  \textrm{metric}
   \hskip-4mm\end{array}
 \right|
\overset{\xi\,\textrm{-Killing}}\longrightarrow
  \left|
   \begin{array}{c}\hskip-2mm
  \textrm{weak} \\
  \textrm{K-contact}
   \hskip-1mm\end{array}
 \right|
\overset{{\cal N}^{\,(1)}=0}
\longrightarrow
 \left|
   \begin{array}{c}\hskip-2mm
 \textrm{weak} \\
 \textrm{Sasakian}
   \hskip-1mm\end{array}
 \right|.
\]
In \cite{RWo-2}, we prove that in the case of a weak structure, the partial Ricci flow,
$\partial g/\partial t = -2\,(\Ric^\bot(g)-\Phi\,{\rm id}^\bot)$
where $\Phi:M\to\RR$,
reduces to ODE's,
$(d/dt)\Ric^\bot = 4\Ric^\bot(\Ric^\bot -\Phi\,{\rm id}^\bot)$,
and that weak structures with positive partial Ricci curvature retract to the classical structures with positive constant partial Ricci~curvature.

We define a {weak almost Hermitian structure} on a Riemannian manifold of even dimension, $(M^{2n}, g)$, equipped with a skew-symmetric {\rm (1,1)}-tensor $f$
by condition
$f^{\,2}<0$.
Such~$(g, f)$ is a {weak K\"{a}hlerian structure}, if $\nabla f=0$, where $\nabla$ is the Levi-Civita connection,
and a {weak nearly K\"{a}hlerian structure}, if $(\nabla_X f)X=0$.
Weak nearly Sasakian/cosymplectic manifolds are defined by a similar condition \eqref{E-nS-01} together with \eqref{2.2}
in the same spirit starting from Sasakian and cosymplectic manifolds, see~\cite{rov-122}:

A distribution ${\cal D}\subset TM$ is \textit{totally geodesic} if and only if
$\nabla_X Y+\nabla_Y X\in{\cal D}\ (X,Y\in{\cal D})$,
this is the case when {any geodesic of $M$ that is tangent to ${\cal D}$ at one point is tangent to ${\cal D}$ at all its points},
e.g., \cite{R2021}.
Any such integrable distribution defines a totally geodesic foliation.
A~foliation, whose orthogonal distribution is totally geodesic, is said to be a Riemannian~foliation.

At first glance, the results in \cite{RWo-2,rov-119} suggest that weak structures are not ``far'' from classical ones.
But it turns out that weak nearly Sasakian/cosymplectic manifolds are relatively ``far'' from classical~ones.
In~\cite{rov-122,rov-129,rov-128}, we gave examples of $(4n+1)$-dimensional proper weak nearly Sasakian mani\-folds
and found conditions:
\begin{eqnarray}
\label{E-nS-10}
 && (\nabla_X\,\widetilde Q)\,Y=0\quad (X\in TM,\ Y\in\ker\eta),\\
\label{E-nS-04c}
 && R_{\widetilde Q X,Y}Z\in\ker\eta\quad (X,Y,Z\in\ker\eta),
\end{eqnarray}
(trivial for~$\widetilde Q=0$) under which weak nearly cosymplectic manifolds split and weak nearly Sasakian manifolds are Sasakian
-- which generalizes the results in~\cite{NDY-2018}.
We~do not extend \eqref{E-nS-10} for $Y=\xi$, since then at any point either $\nabla\xi=0$ or~$\widetilde Q=0$.
By \eqref{E-nS-04c} and the first Bianchi identity, we obtain $R_{\,X,Y}\,\widetilde Q Z\in\ker\eta\ (X,Y,Z\in\ker\eta)$.

The notion of a warped product is popular in differential geometry as well as in general rela\-tivity.
Some solutions of Einstein field equations are warped products and some spacetime models
are warped products.
Z.~Olszak \cite{olszak1991normal} characterized in terms of warped products a class of almost contact
manifolds, known as $\beta$-Kenmotsu manifolds for $\beta=const>0$, given~by
\begin{equation}\label{2.3-patra}
 (\nabla_{X} f)Y=\beta\{g(f X, Y)\,\xi -\eta(Y)\,f X\},
\end{equation}
and defined by K. Kenmotsu \cite{kenmotsu1972class}, when $\beta=1$.
In \cite{rov-126} we show that a weak $\beta$-Kenmotsu manifold is locally the warped product $(-\varepsilon,\varepsilon)\times_\sigma\bar M$, where $(\partial_t\,\sigma)/\sigma=\beta\ne0$, and $\bar M$ has a weak Hermitian structure,
and extend some results on Ricci-type solitons and Einstein-type metrics from Kenmotsu
to weak $\beta$-Kenmotsu~case.

Several authors studi\-ed the problem of finding skew-symmetric parallel 2-tensors (different from almost complex structures)
on a Riemannian manifold and classified them, e.g., \cite{H-2022},
 or proved that some spaces do not admit them, e.g.,~\cite{KKN}.
The~idea of considering the bundle of almost-complex structures compatible with a given metric led to the twistor construction
and then to twistor string theory.
Thus, we delegate to further work the study of the bundles of weak Hermitian and weak almost contact structures that are compatible with a given metric.

\smallskip

The article surveys recent results in \cite{RWo-2,rov-119,rov-126,rov-117,rov-122,rov-128,rov-129,rov-130,r-2024}.

In Section~\ref{sec:02}, following the introductory Section~\ref{sec:01}, we recall some results regarding weak almost-contact mani\-folds.
In Section~\ref{sec:03}, we present rigidity~results for weak Sasakian manifolds
and characterize weak cosymplectic manifolds in the class of weak almost contact metric manifolds.
In Section~\ref{sec:04}, we generalize some properties of $K$-contact manifolds to weak $K$-contact case
and characterize weak $K$-contact manifolds among all weak contact metric manifolds.
In~Section~\ref{sec:05}, we study geometry of weak nearly Sasakian/co\-symplectic structures and
weak nearly Sasakian/cosymplectic
hypersurfaces in weak nearly K\"{a}hler manifolds (generalizing nearly K\"{a}hler manifolds).
In Section~\ref{sec:06}, we show that a weak nearly cosymplectic manifold with parallel Reeb vector field is locally the Riemannian product of a real line and a weak nearly K\"{a}hler~manifold.
In~Section~\ref{sec:07}, we present examples of proper weak nearly Sasakian manifolds,
show that weak nearly Sasakian manifolds satisfying conditions \eqref{E-nS-10}--\eqref{E-nS-04c} have a foliated~structure and
give two theorems characterizing Sasakian manifolds in the class of weak almost contact metric manifolds.
In Section~\ref{sec:08}, we study weak $\beta$-Kenmotsu manifolds.

\section{Preliminaries}
\label{sec:02}
 \vglue-10pt
 \indent

Here, we discuss the basic properties of weak structures.

By
\eqref{E-nS-2.1}, $\ker\eta$ is $Q$-invariant
and the following equalities are true:
\begin{eqnarray*}
 {f}\,\xi=0,\quad \eta\circ{f}=0,\quad \eta\circ Q=\eta,\quad [Q,\,{f}]:={Q}\circ{f} - {f}\circ{Q}=0.
\end{eqnarray*}
Recall the following formulas with the Lie derivative $\pounds$ in the $Z$-direction:
\begin{eqnarray*}
 && (\pounds_{Z}{f})X = [Z, {f} X] - {f} [Z, X],\quad
 (\pounds_{Z}\,\eta)X = Z(\eta(X)) - \eta([Z, X]) , \\
 && (\pounds_{Z}\,g)(X,Y)
 = g(\nabla_{X}\,Z, Y) + g(\nabla_{Y}\,Z, X).
\end{eqnarray*}
The Nijenhuis torsion
of ${f}$ and the exterior derivative
of $\eta$ and $\Phi$ are given~by
\begin{eqnarray*}
 && [{f},{f}](X,Y) = {f}^2 [X,Y] + [{f} X, {f} Y] - {f}[{f} X,Y] - {f}[X,{f} Y],
 \\
 && d\eta(X,Y) = (1/2)\,\{X(\eta(Y)) - Y(\eta(X)) - \eta([X,Y])\} ,\\
 && d\Phi(X,Y,Z) = ({1}/{3})\,\big\{ X\,\Phi(Y,Z) + Y\,\Phi(Z,X) + Z\,\Phi(X,Y) \notag\\
 &&\qquad -\,\Phi([X,Y],Z) - \Phi([Z,X],Y) - \Phi([Y,Z],X)\big\}.
\end{eqnarray*}
The following tensors on almost contact manifolds are well known, see~\cite{blair2010riemannian}:
\begin{eqnarray*}
 && {\cal N}^{\,(2)}(X,Y) = (\pounds_{{f} X}\,\eta)Y - (\pounds_{{f} Y}\,\eta)X
 = 2\,d\eta({f} X,Y) - 2\,d\eta({f} Y,X),  \\
 && {\cal N}^{\,(3)}(X) = (\pounds_{\xi}{f})X
 = [\xi, {f} X] - {f} [\xi, X],\\
 && {\cal N}^{\,(4)}(X) = (\pounds_{\xi}\,\eta)X
 = \xi(\eta(X)) - \eta([\xi, X]) = 2\,d\eta(\xi, X).
\end{eqnarray*}

The following statement is based on \cite[Theorem~6.1]{blair2010riemannian}, i.e., $Q={\rm id}_{\,TM}$.

\begin{proposition}
Let a weak almost contact structure be normal. Then ${\cal N}^{\,(3)}={\cal N}^{\,(4)}=0$ and
 ${\cal N}^{\,(2)}(X,Y) =\eta([\widetilde QX,\,{f} Y])$;
moreover, the vector field $\xi$ is geodesic.
\end{proposition}

\begin{proposition}
For a metric weak almost contact structure
we get
\begin{eqnarray*}
 && 2\,g((\nabla_{X}{f})Y,Z) = 3\,d\Phi(X,{f} Y,{f} Z) - 3\, d\Phi(X,Y,Z) + g({\cal N}^{\,(1)}(Y,Z),{f} X)\notag\\
 && \quad +\,{\cal N}^{\,(2)}(Y,Z)\,\eta(X) + 2\,d\eta({f} Y,X)\,\eta(Z) - 2\,d\eta({f} Z,X)\,\eta(Y)
 + {\cal N}^{\,(5)}(X,Y,Z),
\end{eqnarray*}
where a skew-symmet\-ric with respect to $Y$ and $Z$ (0,3)-tensor ${\cal N}^{\,(5)}$ is defined by
\begin{eqnarray}\label{E-N5}
\nonumber
 && {\cal N}^{\,(5)}(X,Y,Z) = ({f} Z)\,(g(X, \widetilde QY)) -({f} Y)\,(g(X, \widetilde QZ)) +g([X, {f} Z], \widetilde QY)\\
 &&\qquad -\,g([X,{f} Y], \widetilde QZ) +g([Y,{f} Z] -[Z, {f} Y] - {f}[Y,Z],\ \widetilde Q X).
\end{eqnarray}
For particular values of the tensor ${\cal N}^{\,(5)}$ we get
\begin{eqnarray*}
\nonumber
 && {\cal N}^{\,(5)}(X,\xi,Z) = -{\cal N}^{\,(5)}(X, Z, \xi) = g( {\cal N}^{\,(3)}(Z),\, \widetilde Q X),\\
\nonumber
 && {\cal N}^{\,(5)}(\xi,Y,Z) = g([\xi, {f} Z], \widetilde QY) -g([\xi,{f} Y], \widetilde QZ),\\
 && {\cal N}^{\,(5)}(\xi,\xi,Y) = {\cal N}^{\,(5)}(\xi,Y,\xi)=0.
\end{eqnarray*}
\end{proposition}


\begin{theorem}
{\rm (i)}~On a weak ${K}$-contact manifold the vector field $\xi$ is Killing and geodesic, and
\begin{eqnarray*}
 && {\cal N}^{\,(1)}(\xi,\,\cdot)={\cal N}^{\,(5)}(\xi,\,\cdot\,,\,\cdot)={\cal N}^{\,(5)}(\,\cdot\,,\xi,\,\cdot)={\cal N}^{\,(5)}(\cdot\,, \xi, \,\cdot)=0,\\
 && \pounds_{\,\xi}{Q}=\nabla_\xi{Q}=0,\ \nabla_{\xi}{f}=0,\quad \nabla_X\,\xi=-f X ,\\
 && g((\nabla_{X}{f})Y,Z) = d\eta({f} Y,X)\,\eta(Z) - d\eta({f} Z,X)\,\eta(Y)
 +\frac12\,\eta([\widetilde Q Y,\,{f} Z])\,\eta(X) + \frac12\,{\cal N}^{\,(5)}(X,Y,Z).
\end{eqnarray*}

{\rm (ii)}~For a weak almost Sasakian structure, $\xi$ is a geodesic vector field and the tensors ${\cal N}^{\,(2)}$ and ${\cal N}^{\,(4)}$ vanish;
moreover, ${\cal N}^{\,(3)}$ vanishes if and only if $\,\xi$ is a Killing vector field, and
\begin{eqnarray}\label{3.1A}
 g((\nabla_{X}{f})Y,Z) = \frac12\,g({\cal N}^{\,(1)}(Y,Z),{f} X)
 +g(fX,fY)\eta(Z) -g(fX,fZ)\eta(Y) {+} \frac12\,{\cal N}^{\,(5)}(X,Y,Z).
\end{eqnarray}
In particular, $g((\nabla_{\xi}{f})Y,Z) = \frac12\,{\cal N}^{\,(5)}(\xi,Y,Z)$.

{\rm (iii)}~For a weak almost cosymplectic manifold, we get ${\cal N}^{\,(2)}={\cal N}^{\,(4)}=0$, ${\cal N}^{\,(1)}=[{f},{f}]$, and
$\xi$ is geodesic. Moreover, ${\cal N}^{\,(3)}=0$ if and only if $\,\xi$ is a Killing vector~field.
\end{theorem}

\section{The rigidity of Sasakian structure \\ and characteristic of a cosymplectic structure}
\label{sec:03}
 \vglue-10pt
 \indent

Here, we present rigidity~results for weak Sasakian manifolds
and characterize weak cosymplectic manifolds in the class of weak almost contact metric manifolds.


\begin{proposition}
For a weak Sasakian structure, \eqref{3.1A} reduces to
\begin{eqnarray*}
 && g((\nabla_{X}{f})Y,Z) = g(QX,Y)\,\eta(Z) - g(QX,Z)\,\eta(Y) \notag\\
 && +\,\eta(X)\big(\eta(Y)\,\eta(Z) - \eta(Y)\,\eta(Z)\big)  +\frac{1}{2}\,{\cal N}^{\,(5)}(X,Y,Z) .
\end{eqnarray*}
\end{proposition}

The equality ${\cal N}^{(3)}=0$ is valid for weak Sasakian manifolds, since it is true for Sasakian manifolds, see Theorem~\ref{T-4.1}.
We get the rigidity of the Sasakian~structure.

\begin{theorem}\label{T-4.1}
A weak almost contact metric structure is a weak Sasakian structure if and only if it is a Sasakian structure.
\end{theorem}


\begin{proposition}
Let $({f},Q,\xi,\eta,g)$ be a weak
cosymplectic structure. Then
\begin{eqnarray}\label{6.1}
 2\,g((\nabla_{X}{f})Y,Z) = {\cal N}^{\,(5)}(X,Y,Z).
\end{eqnarray}
In particular, using \eqref{6.1}
 and \eqref{E-nS-2.1}, we get
 $g(\nabla_{X}\,\xi,\,Q Z) = -\frac12\,{\cal N}^{\,(5)}(X,\xi,{f} Z)$.
\end{proposition}

Recall that a ${K}$-structure is a cosymplectic structure if and only if ${f}$ is parallel.
 The following our theorem generalizes this result.

\begin{theorem}\label{thm6.2D}
A weak almost contact metric structure with $\nabla{f}=0$
is a~weak cosymplectic structure with ${\cal N}^{\,(5)}=0$.
\end{theorem}

\begin{example}\label{Ex-3.1}\rm
Let $M$ be a $2n$-dimensional smooth manifold and $\bar{f}:TM\to TM$ an endomorphism of rank $2n$ such that $\nabla\bar{f}=0$.
To construct a weak cosymplectic structure on $M\times\mathbb{R}$ or $M\times \mathbb{S}^1$,
take any point $(x, t)$ of either space and set $\xi = (0, d/dt)$, $\eta =(0, dt)$~and
\[
 {f}(X, Y) = (\bar{f} X, 0),\quad
 Q(X, Y) = (-\bar{f}^{\,2} X,\, Y).
\]
where $X\in M_x$ and $Y\in\{\mathbb{R}_t, \mathbb{S}^1_t\}$. Then \eqref{E-nS-2.1} holds and Theorem~\ref{thm6.2D} can be used.
\end{example}


\section{Weak $K$-contact structure}
\label{sec:04}
 \vglue-10pt
 \indent

Here, we generalize some properties of $K$-contact manifolds to weak $K$-contact case.
We characterize weak $K$-contact manifolds among all weak contact metric mani\-folds
by the following well known property of $K$-contact manifolds, see \cite{blair2010riemannian}:
\begin{equation}\label{E-30}
 \nabla\,\xi = -{f} .
\end{equation}

\begin{theorem}
A weak contact metric manifold is weak $K$-contact $($that is $\xi$ is a Killing vector field$)$ if and only if \eqref{E-30} is valid.
\end{theorem}


A~Riemannian manifold with a unit Killing vector field and the property
$R_{{X},\,{\xi}}\,{\xi}
={X}\ ({X}\bot\,\xi)$ is a $K$-contact manifold, e.g.,
\cite{blair2010riemannian}.
We~generalize this
in the following

\begin{theorem}\label{prop2.1b}
A Riemannian manifold $(M^{\,2n+1},g)$ admitting a unit Killing field $\xi$ with positive
$\xi$-sectional curvature
is a weak $K$-contact manifold $M({f},Q,\xi,\eta,g)$ with
the tensors: $\eta = g(\cdot,\, \xi)$,
${f} = -\nabla\,\xi$, see \eqref{E-30}, and $Q {X} = R_{{X},\,{\xi}}\,{\xi}\ ({X}\in\ker\eta)$.
\end{theorem}

\begin{example}\rm
By Theorem~\ref{prop2.1b}, we can search for examples of weak $K$-contact (not $K$-contact) manifolds
among Riemannian manifolds of positive sectional curvature that admit unit Killing vector fields.
Indeed, let $M$ be a convex hypersurface (ellipsoid) with induced metric $g$ of the Euclidean space~$\mathbb{R}^{2n+2}$,
\[
 M = \Big\{(u_1,\ldots,u_{2n+2})\in\mathbb{R}^{2n+2}: \sum\nolimits_{\,i=1}^{n+1} u_i^2 + a\sum\nolimits_{\,i=n+2}^{2n+1} u_i^2 = 1\Big\},
 \quad 0<a=const\ne1,
\]
where $n\ge1$ is odd. The sectional curvature of $(M,g)$ is positive. It follows that
\[
 \xi = (-u_2, u_1, \ldots , -u_{n+1}, u_{n}, -\sqrt a\,u_{n+3}, \sqrt a\,u_{n+2}, \ldots , -\sqrt a\,u_{2\,n+2}, \sqrt a\,u_{2\,n+1})
\]
is a Killing vector field on $\mathbb{R}^{2\,n+2}$, whose restriction to $M$ has unit length.
Since $\xi$ is tangent to
$M$,
so $\xi$ is a unit Killing vector field on $(M,g)$, see \cite[p.~5]{D-B-2021}.
For $n=1$, we get a weak $K$-contact manifold
$M^3=\big\{u_1^2 + u_2^2 + a u_3^2 + a u_4^2 = 1\big\}\subset\mathbb{R}^{4}$ with $\xi = (-u_2, u_1, -\sqrt a\,u_{4}, \sqrt a\,u_{3})$.
\end{example}

\begin{proposition}\label{Pr6.2D}
For a weak $K$-contact manifold, the following equalities hold:
\begin{eqnarray*}
 R_{{\xi},{X}} = \nabla_{X}{f}, \quad
 R_{\,{\xi},{X}}\,{\xi} = {f}^2 {X} ,\quad
 {\rm Ric}({\xi},{\xi}) = \tr Q = 2\,n + \tr\tilde Q .
\end{eqnarray*}
\end{proposition}

If a Riemannian manifold
admits a unit Killing vector field $\xi$, then
$K(\xi,{X})\ge0\ ({X}\,\perp\,\xi,\ {X}\ne0)$, thus ${\rm Ric}({\xi},{\xi})\ge0$; moreover, ${\rm Ric}({\xi},{\xi})\equiv0$ if and only if
$\nabla\xi\equiv0$.
In the case of $K$-contact manifolds, $K(\xi,{X})=1$, see \cite[Theorem~7.2]{blair2010riemannian}.

\begin{corollary}\label{Cor_K-positive}
For a weak $K$-contact manifold, the $\xi$-sectional curvature is
\begin{equation*}
 K(\xi,{X})=g(Q{X},{X})>0\quad ({X}\in{\cal D},\ \|{X}\|=1);
\end{equation*}
therefore, for the Ricci curvature we get ${\rm Ric}({\xi},{\xi})>0$.
\end{corollary}

Using Theorem~\ref{prop2.1b} and Corollary~\ref{Cor_K-positive},
we show that
a weak $K$-contact structure
can be deformed (by the partial Ricci flow) to a $K$-contact structure, see~\cite{RWo-2}.

\begin{corollary}
A~weak $K$-contact manifold $M({f},Q,\xi,\eta,g_0)$ admits a smooth fa\-mily of metrics $g_t\ (t\in\mathbb{R})$,
such that $M({f}_t,Q_t,\xi,\eta,g_t)$ are weak $K$-contact manifolds with certainly defined ${f}_t$ and $Q_t$;
moreover, $g_t$ converges exponentially fast, as $t\to-\infty$, to a limit metric $\hat g$ that gives a $K$-contact structure.
\end{corollary}

The following theorem generalizes a well known result, e.g., \cite[Proposition 5.1]{YK-1985}.

\begin{theorem}
A weak $K$-contact manifold with conditions $(\nabla\Ric)(\xi,\,\cdot)=0$
and $\tr Q=const$ is an Einstein manifold of scalar curvature $r=(2\,n + 1)\tr Q$.
\end{theorem}

\begin{remark}\rm
For a weak $K$-contact manifold, by $\nabla_{\xi}{f}=\frac12\,{\cal N}^{\,(5)}(\xi,{Y},\,Z)=0$ and Proposition~\ref{Pr6.2D},
we get the equality (well known for $K$-contact manifolds, e.g., \cite{blair2010riemannian}):
 $\Ric^\sharp(\xi)=\sum\nolimits_{\,i=1}^{\,2n}(\nabla_{e_i}{f})\,e_i$,
where $(e_i)$ is any local orthonormal basis of~$\ker\eta$;
and for contact manifolds we have $\sum\nolimits_{\,i=1}^{\,2n}(\nabla_{e_i}{f})\,e_i=2\,n\,\xi$.
For $K$-contact manifolds, this gives $\Ric^\sharp(\xi)=2\,n\,\xi$,
and $\Ric(\xi,\xi)=2\,n$; moreover, the last condition characterizes $K$-contact mani\-folds among all contact metric manifolds.
\end{remark}

The Ricci curvature of any $K$-contact manifold satisfies the condition
\begin{equation}\label{E-K-Ric-{X}}
 \Ric(\xi,{X}) = 0\quad ({X}\in{\cal D}).
\end{equation}

Quasi Einstein manifolds are defined
by the condition
 $\Ric({X}, {Y}) = a\,g({X}, {Y}) + b\,\mu({X})\,\mu({Y})$,
where $a$ and $b\ne0$ are real scalars, and $\mu$ is a 1-form of unit norm.

The next our theorem generalizes \cite[Theorem 3.1]{G-D-2020}.

\begin{theorem}
Let a weak $K$-contact manifold with $\tr Q=const$ satisfy
\eqref{E-gg-r-e2}
with $c_1 a\ne -1$ for $a=\lambda+c_2\tr Q$.
If~\eqref{E-K-Ric-{X}} is true, then $f=const$. Furthermore,

\noindent\
-- if $c_1 a\ne0$, then
 ${\rm Hess}_{f_2} = \frac1a\,df_2\otimes df_2 - \frac{c_2}{c_1 a}\,\Ric - \frac\lambda{c_1 a}\,g$;
if $c_1 a\ne-1$, then $f_2=const$; and if $c_2\ne0$, then $(M,g)$ is an Einstein manifold.

\noindent\
-- if $a=0$ and $c_1\ne0$, then
 $0 = c_2\Ric -c_1\,df_2\otimes df_2 + \lambda\,g$.
If~$c_2\ne0$ and $f_2\ne const$, then we get a gradient quasi Einstein manifold.

\noindent\
-- for $c_1=0$, then
$c_2\Ric + \lambda\,g=0$, and for $c_2\ne0$ we get an Einstein~manifold.
\end{theorem}

Quasi contact metric manifolds (introduced by Y.~Tashiro)
are an extension of contact metric manifolds.
In \cite{r-2024}, we study a weak analogue of quasi contact metric manifolds and provide new criterions for K-contact and Sasakian manifolds.

We pose the following {questions}.
{Is the condition ``the $\xi$-sectional curvature is positive" sufficient for a weak almost contact metric manifold to be weak K-contact}?
{Does a weak almost contact metric manifold of dimension $>3$ have some positive $\xi$-sectional curvature}?
{Is a compact weak K-contact Einstein manifold a Sasakian manifold}?
{When a~weak almost contact manifold equipped with a Ricci-type soliton structure,
carries a canonical (for example, with constant sectional curvature or Einstein-type) metric}?
Is a compact {weak} K-contact Einstein mani\-fold a~Sasakian manifold?
Is a compact {weak} K-contact mani\-fold admitting a generalized Ricci soliton structure a~Sasakian manifold?
To answer these questions, we need to generalize some deep results about contact manifolds to weak contact~manifolds.


\section{Weak Nearly Sasakian/Cosymplectic Manifolds}
\label{sec:05}
 \vglue-10pt
 \indent

Here, we study geometry of weak nearly Sasakian/cosymplectic manifolds
and hypersurfaces in weak nearly K\"{a}hler manifolds (generalizing nearly K\"{a}hler manifolds)

The following result generalizes {Proposition~3.1} in~\cite{blair1976} and {Theorem~5.2} in~\cite{blair1974}.

\begin{theorem}
(i)~Both on weak nearly Sasakian and weak nearly cosymplectic manifolds the vector field $\xi$ is geodesic;
moreover, if the condition \eqref{E-nS-10} is valid, then the vector field $\xi$ is Killing.
(ii)~There are no weak nearly cosymplectic structures with the condition \eqref{E-nS-10}, which are weak contact metric structures.
\end{theorem}

\begin{theorem}
For a weak nearly cosymplectic manifold, $\nabla\xi=0$ if and only if the manifold
is locally
a metric
product of $\RR$ and a weak nearly K\"{a}hler~manifold.
\end{theorem}

\begin{proposition}\label{P-22}
Let a weak almost contact manifold $M({f},Q,\xi,\eta)$ satisfy \eqref{E-nS-10} and
 $Q\,|_{\,\ker\eta}=\lambda\,{\rm id}|_{\,\ker\eta}$
{for} a positive function $\lambda\in C^\infty(M)$.
Then, $\lambda=const$ and
\begin{equation}\label{Tran'}
 {f} = \sqrt\lambda\,\tilde{f} .
\end{equation}
Moreover, if a weak nearly Sasakian/cosymplectic
structure $({f},Q,\xi,\eta,g)$
satisfies \eqref{Tran'} and
\begin{equation}\label{Tran2'}
 g|_{\,\ker\eta} = \lambda^{\,-\frac12}\,\tilde g|_{\,\ker\eta},\quad
 g(\xi,\,\cdot) = \tilde{g}(\xi,\,\cdot) ,
\end{equation}
{then} $(\tilde{f},\xi,\eta, \tilde{g})$ is a nearly Sasakian/cosymplectic
structure.
\end{proposition}

\begin{example}\rm
Let $M({f},Q,\xi,\eta,g)$ be a three-dimensional weak almost-contact metric ma\-nifold.
The tensor $Q$ has on the plane field $\ker\eta$ in the form $\lambda\,{\rm id}_{\,\ker\eta}$ for some positive function $\lambda\in C^\infty(M)$.
Suppose that
\eqref{E-nS-10} is true, then $\lambda=const$ and this structure reduces to the almost-contact metric structure $(\tilde{f},\xi,\eta,\tilde g)$ satisfying \eqref{Tran'} and \eqref{Tran2'}.

Let \eqref{E-nS-01} hold for $M({f},Q,\xi,\eta,g)$.
By Proposition~\ref{P-22}(ii), $M(\tilde{f},\xi,\eta,\tilde g)$ is nearly cosymplectic or nearly Sasakian, respectively.
Since $\dim M=3$, we obtain Sasakian
({Theorem~5.1} in~\cite{Ol-1980})
or cosymplectic
(see~\cite{JKK-94})
structures $(\tilde{f},\xi,\eta,\tilde g)$, respectively.
\end{example}

We generalize rigidity {Theorem~3.2} in~\cite{blair1976}.

\begin{theorem}
For a weak nearly Sasakian structure satisfying \eqref{E-nS-10}, normality $({N}^{\,(1)}=0)$ is equivalent to a weak contact metric property $(d\eta=\Phi)$.
Therefore, a normal weak nearly Sasakian structure satisfying \eqref{E-nS-10} is Sasakian.
\end{theorem}

\begin{proposition}
A 3-dimensional weak nearly cosymplectic structure
satisfy\-ing \eqref{E-nS-10} reduces to cosymplectic one.
\end{proposition}

\begin{example}\rm
Let a 3-dimensional weak nearly Sasakian manifold $M({f},Q,\xi,\eta,g)$ satis\-fy \eqref{E-nS-10}.
By \eqref{E-nS-2.1},
$Q$ has on the plane field $\ker\eta$ the form $\lambda\,{\rm id}_{\,\ker\eta}$ for some positive $\lambda\in\mathbb{R}$.
This structure reduces to the nearly Sasakian structure $(\tilde{f},\xi,\eta,\tilde g)$, where
 $\tilde{f} = \lambda^{\,-\frac12}\,{f}$,
 $\tilde  g|_{\,\ker\eta} = \lambda^{\,\frac12}\,g|_{\,\ker\eta}$,
 $\tilde g(\xi,\,\cdot) = {g}(\xi,\,\cdot)$.
Since $\dim M=3$, the structure $(\tilde{f},\xi,\eta,\tilde g)$ is Sasakian.
\end{example}

Next, we will study weak nearly Sasakian/cosymplectic hypersurfaces in
weak nearly K\"{a}hler manifolds (generalizing nearly K\"{a}hler manifolds).

\begin{example}
\rm
Let $(\bar M, \bar{f}, \bar g)$ be a weak nearly K\"{a}hler manifold:$(\bar\nabla_X\bar{f})X=0\ (X\in T\bar M)$.
To build a weak nearly cosymplectic structure $({f}, Q,\xi,\eta,g)$ on the product
$M=\bar M\times\mathbb{R}$
of $(\bar M, \bar g)$ and a Euclidean line $(\mathbb{R},\partial_t)$, we take any point $(x, t)$ of $M$ and set
\[
 \xi = (0, \partial_t),\quad
 \eta =(0, dt),\quad
 {f}(X, \partial_t) = (\bar{f} X, 0),\quad
 Q(X, \partial_t) = (-\bar{f}^{\,2} X, \partial_t),\quad X\in T_x\bar M
\]
(similarly to Example~\ref{Ex-3.1}).
Note that if $\bar\nabla_X\bar{f}^2=0\ (X\in T\bar M)$, then \eqref{E-nS-10} holds.
\end{example}

The scalar second fundamental form $b$ of a
hypersurface $M\subset (\bar M, \bar g)$ with a unit normal $N$ is related with $\overline\nabla$ and the Levi-Civita connection $\nabla$ induced on the $M$ metric $g$ via the Gauss equation
\begin{equation*}
 \overline\nabla_X Y = \nabla_X Y +b(X,Y)\,N\quad (X,Y\in TM).
\end{equation*}
{The} Weingarten operator $A_N: X\mapsto -\overline\nabla_X N$ is related with $b$ via the equality
 $\bar g(b(X,Y),N)=g(A_N(X),Y)\ (X,Y\in TM)$.
{A} hypersurface is called \textit{totally geodesic} if $b=0$.
A hypersurface is called \textit{quasi-umbilical}~if
 $b(X,Y) = c_1\,g(X,Y)+ c_2\,\mu(X)\,\mu(Y)$,
where $c_1,c_2\in C^\infty(M)$ and $\mu\ne0$ is a one-form.

\begin{lemma}
A hypersurface $(M,g)$ with a unit normal $N$
in a weak Hermitian manifold $(\bar M, \bar{f}, \bar g)$
inherits a weak almost-contact structure $({f},Q,\xi,\eta,g)$ given by
\[
 \xi= \bar{f}\,N,\quad
 \eta = \bar g(\bar{f}\,N, \,\cdot),\quad
 {f} = \bar{f} + \bar g(\bar{f}\,N, \,\cdot)\,N,\quad
 Q = -\bar{f}^{\,2}  + \bar g(\bar{f}^{\,2} N,\,\cdot)\,N.
\]
\end{lemma}

The following theorem generalizes the fact (see~\cite{blair1976}) that a hypersurface of a nearly K\"{a}hler manifold is nearly Sasakian or nearly cosymplectic if and only if it is quasi-umbilical with respect to the (almost) contact form.

\begin{theorem}
Let $M$ be a hypersurface with a unit normal $N$ of a weak nearly K\"{a}hler manifold $(\bar M^{2n+2}, \bar{f}, \bar g)$.
Then, the induced structure $({f},Q,\xi,\eta,g)$ on $M$ {is} 
\begin{equation*}
(i)~\mbox{Weak nearly Sasakian}; \quad
(ii)~\mbox{Weak nearly cosymplectic}.
\end{equation*}
{This is true if} 
and only if $M$ is quasi-umbilical with the
scalar 2nd fundamental~form
\begin{equation*}
 (i) \ b(X,Y) = g(X,Y)+ (b(\xi,\xi)-1)\,\eta(X)\,\eta(Y);\ \
 (ii)\ b(X,Y) = b(\xi,\xi)\,\eta(X)\,\eta(Y).
\end{equation*}
{In} both cases, $A_N{f}+{f} A_N=2{f}$ is true, and \eqref{E-nS-10} follows from the condition
\[
 ((\overline\nabla_X\,\bar{f}^2)Y)^\top=0\quad (X,Y\in TM,\ Y\perp\xi).
\]
\end{theorem}

\section{Splitting of weak nearly cosymplectic manifolds}
\label{sec:06}

 \vglue-10pt
 \indent

Here, we show
that a weak nearly cosymplectic manifold satisfies \eqref{E-nS-04cc}, if we assume a weaker condition \eqref{E-nS-04c}.
Then, we present the splitting of weak nearly cosymplectic manifolds satisfying \eqref{E-nS-10} and~\eqref{E-nS-04c}
and generalize some well known results.
We also characterize 5-dimensional weak nearly cosymplectic mani\-folds.

We define a (1,1)-tensor $h = \nabla\xi$ on $M$ as in the classical case, e.g., \cite{E-2005}.
Note that $h=0$ if and only if
$\ker\eta$ is integrable, i.e., $[X,Y]\in\ker\eta\ (X,Y\in\ker\eta)$.
Since $\xi$ is a geodesic vector field,
we get $h\,\xi=0$ and $h(\ker\eta)\subset\ker\eta$.
Since $\xi$ is a Killing vector field, the tensor $h$ is skew-symmetric:
 $g(h X,\, X) = g(\nabla_X\,\xi, X) = \frac12\,(\pounds_\xi\,g)(X,X) = 0$.
We~also get $\eta\circ h = 0$ and
  $d\,\eta(X,\,\cdot) = \nabla_X\,\eta = g(hX,\,\cdot)$.


\begin{lemma}
For a weak nearly cosymplectic manifold $M({f},Q,\xi,\eta,g)$ we obtain
\begin{eqnarray*}
 (\nabla_X\,h)\,\xi = -h^2 X,\qquad
 (\nabla_X\,{f})\,\xi = -{f}\,h X .
\end{eqnarray*}
Moreover, if the condition \eqref{E-nS-10} is true, then
\begin{eqnarray}
\label{E-nS-01b}
 && h\,{f} + {f}\,h =0\quad (h\ {\rm anticommutes\ with}\ {f}),\\
\label{E-nS-01d}
 && h\,Q = Q\,h\quad (h\ {\rm commutes\ with}\ Q).
\end{eqnarray}
\end{lemma}

\begin{lemma}
For a weak nearly cosymplectic manifold we get
the equality
\begin{equation*}
 g(R_{{f} X,Y}Z, V) +g(R_{X,{f} Y}Z, V) +g(R_{X,Y}{f} Z, V) +g(R_{X,Y}Z, {f} V) = 0.
\end{equation*}
Moreover, if the conditions \eqref{E-nS-10} and \eqref{E-nS-04c} are true, then
 $g(R_{\,\xi, Z}\,{f} X,{f} Y) = 0$.
\end{lemma}

\begin{lemma}
Let a weak nearly cosymplectic manifold satisfy \eqref{E-nS-10}--\eqref{E-nS-04c}, then
\begin{equation*}
 g((\nabla_{X}\,{f})Y, {f} h Z) = \eta(X)\,g(hY, hQZ) -\eta(Y)\,g(hX, hQZ).
\end{equation*}
\end{lemma}

\begin{lemma}
For a weak nearly cosymplectic manifold satisfying \eqref{E-nS-04c}, we get
\begin{eqnarray}\label{E-3.24}
 && (\nabla_X\,h)Y = g(h^2 X, Y)\,\xi - \eta(Y)\,h^2 X , \\
\label{E-3.23}
 && R_{\,\xi, X}Y = -(\nabla_X\,h)Y, \quad
 {\rm Ric}\,(\xi, Z) = -\eta(Z)\,{\rm tr}\,h^2 .
\end{eqnarray}
In particular, $\nabla_\xi\,h=0$ and ${\rm tr}(h^2) = const$.
By \eqref{E-3.24}--\eqref{E-3.23}, we get
\begin{equation*}
  g(R_{\,\xi, X}Y, Z) = -g((\nabla_X\,h)Y, Z) = \eta(Y)\,g(h^2 X, Z) -\eta(Z)\,g(h^2 X, Y).
\end{equation*}
\end{lemma}

The following proposition generalizes \cite[Proposition~4.2]{NDY-2018}.

\begin{proposition}\label{Prop-4.1}
For a weak nearly cosymplectic manifold satisfying \eqref{E-nS-10} and \eqref{E-nS-04c},
the eigenvalues and their multiplicities of the
operator $h^2$ are constant.
\end{proposition}

By Proposition~\ref{Prop-4.1}, the spectrum of the self-adjoint operator $h^2$ has the~form
\begin{equation}\label{E-nS-11b}
 Spec(h^2) = \{0, -\lambda_1^2,\ldots -\lambda_r^2\} ,
\end{equation}
where $\lambda_i$ is a positive real number and $\lambda_i\ne \lambda_j$ for $i\ne j$.
If $X\ne0$ is an eigenvector of $h^2$ with eigenvalue $-\lambda^2_i$, then $X, {f} X, hX$ and $h\,{f} X$ are orthogonal
nonzero eigenvectors of $h^2$ with eigenvalue $-\lambda^2_i$.
Since $h(\xi)=0$, the eigenvalue 0 has multiplicity $2p+1$ for some
$p\ge0$.
Denote by $D_0$ the smooth distribution of the eigenvectors with eigenvalue 0 orthogonal to $\xi$.
Let $D_i$ be the smooth distribution of the eigenvectors with eigenvalue $-\lambda^2_i$.
Thus, $D_0$ and $D_i$ belong to $\ker\eta$ and are ${f}$-invariant and $h$-invariant.
The following proposition generalizes \cite[Proposition~4.3]{NDY-2018}.

\begin{proposition}
Let
a weak nearly cosymplectic manifold satisfy
\eqref{E-nS-10} and \eqref{E-nS-04c},
and let the spectrum of the self-adjoint opera\-tor $h^2$ have the form \eqref{E-nS-11b}. Then,

\noindent\
$(a)$ each distribution $[\xi]\oplus D_i\ (i = 1,\ldots, r)$ is integrable with totally geodesic leaves.

\noindent
Moreover, if the eigenvalue $0$ of $h^2$ is not simple, then

\noindent\
$(b)$ the distribution $D_0$ is integrable with totally geodesic leaves, and each leaf of $D_0$ is endowed with a weak nearly K\"{a}hler structure
$(\bar{f}, \bar g)$
satisfying $\bar\nabla(\bar{f}^{\,2})=0$;

\noindent\
$(c)$ the distribution $[\xi]\oplus D_1\oplus\ldots\oplus D_r$ is integrable with totally geodesic leaves.
\end{proposition}

\begin{proposition}
For a weak nearly cosymplectic (non-weak-cosymplectic) ma\-nifold, $h\equiv0$ if and only if the manifold
is locally isometric to the Riemannian product of a real line and a weak nearly K\"{a}hler (non-weak-K\"{a}hler) manifold.
\end{proposition}

We generalize Theorem~4.5 in \cite{NDY-2018} on splitting of nearly cosymplectic manifolds.

\begin{theorem}
Let $M^{\,2n+1}({f}, Q, \xi, \eta, g)$ be a weak nearly cosymplectic (non-weak-cosymplectic) manifold of dimension $2n+1>5$
with conditions \eqref{E-nS-10} and \eqref{E-nS-04c}.
Then $M$ is locally isometric to one of the Riemannian products:
 $\mathbb{R}\times \bar M^{\,2n}$ or $B^5 \times\bar M^{\,2n-4}$,
where $\bar M(\bar{f},\bar g)$ is a weak nearly K\"{a}hler manifold satisfying $\bar\nabla(\bar{f}^{\,2})=0$,
and $B^5$ is a
weak nearly cosymplectic (non-weak-cosymplectic) manifold satisfy\-ing~\eqref{E-nS-10} and \eqref{E-nS-04c}.
If
$M$ is complete and simply connected, then
the isometry is~global.
\end{theorem}

The following theorem generalizes Theorem~4.4 in \cite{NDY-2018}.

\begin{theorem}
Let $M({f}, Q, \xi, \eta, g)$ be a weak nearly cosymplectic manifold with conditions \eqref{E-nS-10}--\eqref{E-nS-04c}
such that $0$ is a simple eigenvalue of $h^2$. Then $\dim M=5$.
\end{theorem}


\section{Characterization of Sasakian manifolds}
\label{sec:07}
 \vglue-10pt
 \indent

Here, we give examples of proper weak nearly Sasakian manifolds,
present two theorems characterizing Sasakian manifolds in the class of weak almost contact metric manifolds satisfying conditions \eqref{E-nS-10}--\eqref{E-nS-04c}.

On a weak nearly Sasakian manifold satisfying \eqref{E-nS-10}, the unit vector field $\xi$ is Killing ($\pounds_\xi\,g=0$).
Therefore, $\xi$-curves define a Riemannian geodesic foliation.

\begin{example}\rm
We
construct proper weak nearly Sasakian mani\-folds from a pair of classical structures with the same
$\xi,\eta$ and $g$.
Assume $0\le i<n$ and define a
manifold
\[
 M=\{(x_0,x_1,\ldots,x_{4n})\in\mathbb{R}^{4n+1}: x_{4i+2} x_{4i+4}\ne0\}
\]
with standard coordinates $(x_0,x_1,\ldots,x_{4n})$. The vector fields
 $X_0 = \xi = -\partial_0$,
 $X_{4i+1}=2\big(x_{4i+2}\partial_0 - \partial_{4i+1}\big)$,
 $X_{4i+2} = \partial_{4i+2}$,
 $X_{4i+3}=2\big(x_{4i+4}\partial_0 - \partial_{4i+3}\big)$,
 $X_{4i+4} = \partial_{4i+4}$
are pointwise linearly independent.
The non-vanishing Lie brackets are $[X_{4i+1},X_{4i+2}]=[X_{4i+3},X_{4i+4}]=2\,\xi$.
Define a Riemannian metric of $M$ by $g(X_i,X_j)=\delta_{ij}$,
or,
\[
 g = dx_0^2+\sum\nolimits_{\,i}\big\{(1/4-x_{4i+2}^2)dx_{4i+1}^2+dx_{4i+2}^2+(1/4-x_{4i+4}^2)dx_{4i+3}^2+dx_{4i+4}^2\big\}.
\]
Set $\eta=-dx_0$ and define a (1,1)-tensor ${f}_1$ on $M$ by
${f}_1 X_0=0$, ${f}_{1} X_{4i+1}=X_{4i+2}$, ${f}_{1} X_{4i+2}=-X_{4i+1}$, ${f}_{1} X_{4i+3}=X_{4i+4}$,
${f}_1 X_{4i+4}=-X_{4i+3}$.
Thus, $({f}_1,\xi,\eta,g)$ is an almost contact metric structure on $M$.
The non-zero
derivatives $\nabla_{X_{a}}X_{b}$ are
\begin{eqnarray}\label{E-nabla2}
\nonumber
 &\nabla_{X_{4i+1}}X_{4i+2}=-\nabla_{X_{4i+2}}X_{4i+1}=\xi,\quad \nabla_{X_{4i+1}}\xi=\nabla_{\xi}X_{4i+1}=-X_{4i+2}, \\
\nonumber
 & \nabla_{X_{4i+2}}\xi=\nabla_{\xi}X_{4i+2}=X_{4i+1},\quad \nabla_{X_{4i+3}}\xi=\nabla_{\xi}X_{4i+3}=-X_{4i+4},\\
 & \nabla_{X_{4i+4}}\xi=\nabla_{\xi}X_{4i+4}=X_{4i+3},\quad \nabla_{X_{4i+3}}X_{4i+4}=-\nabla_{X_{4i+4}}X_{4i+3}=\xi.
\end{eqnarray}
Thus  $({f}_1,\xi,\eta,g)$ is a Sasakian structure.
In particular, the distribution $\ker\eta$ is curvature invariant.
%
Define a tensor ${f}_2$ on $M$ by
${f}_2 X_{4i+1}=X_{4i+4}$, ${f}_2 X_{4i+4}=-X_{4i+1}$, ${f}_2 X_0=0$, ${f}_2 X_{4i+3}=X_{4i+2}$, ${f}_2 X_{4i+2}=-X_{4i+3}$.
It is easy to check that $({f}_2,\xi,\eta,g)$ is an almost contact metric structure on $M$.
Using \eqref{E-nabla2}, we find that $(\nabla_Y{f}_2)Y=0$ and $(\nabla{f}_2)(X_{4i+1},\xi)=-X_{4i+3} \ne0$, i.e., $({f}_2,\xi,\eta,g)$ is a proper nearly cosymplectic structure.
To construct a weak nearly Sasakian structure using the two above structures, we define a tensor ${f} := \cos(t)\,{f}_1 + \sin(t)\,{f}_2$ for small $t>0$. A tensor $\psi:={f}_1\,{f}_2+{f}_2\,{f}_1$ is self-adjoint with the following nonzero components:
\[
 \psi X_{4i+1} = -2 X_{4i+3},\ \
 \psi X_{4i+2} = -2 X_{4i+4},\ \
 \psi X_{4i+3} = -2 X_{4i+1},\ \
 \psi X_{4i+4} = -2 X_{4i+2}.
\]
Thus $Q={\rm id}-\sin(t)\cos(t)\,\psi\ne {\rm id}$
and
$({f},Q,\tilde\xi,\tilde\eta, \tilde g)$ is a weak nearly Sasakian structure, where $\tilde g=\mu^{2}g$, $\tilde\xi=\mu^{-1}\xi$, $\tilde\eta=\mu\,\eta$ and $\mu=\cos(t)$.
Since the distribution $\ker\eta$ is curvature invariant and $\tilde R=\mu^2 R$,
then $({f},Q,\tilde\xi,\tilde\eta, \tilde g)$ satisfies \eqref{E-nS-04c};
but $({f},Q,\tilde\xi,\tilde\eta, \tilde g)$ does not satisfy \eqref{E-nS-10};
for example, $(\tilde\nabla_{X_{4i+3}}Q)X_{4i+2}
=\sin(2\,t)\,\xi\ne0$.
\end{example}

Here, we generalize some properties of nearly Sasakian manifolds to the case of weak nearly Sasakian manifolds
satisfying \eqref{E-nS-10} and~\eqref{E-nS-04c}.
Define a (1,1)-tensor field $h$ on $M$, as in the classical case, by
 $h = \nabla\xi + {f}$.
We get $\eta\circ h = 0$ and $h(\ker\eta)\subset\ker\eta$. Since $\xi$ is a geodesic field,
 we also get
 $h\,\xi=0$.
Since $\xi$ is a Killing
and ${f}$ is skew-symmetric, the tensor $h$ is skew-symmetric:
\[
 g(h X,\, X) = g(\nabla_X\,\xi, X) + g({f} X, X) = (1/2)\,(\pounds_\xi\,g)(X,X) = 0,
\]
and $\nabla_X\,\eta=g((h-{f})X,\,\cdot)$ holds.
The $\ker\eta$ is integrable
if and only if $h={f}$,
and in this case, our manifold is locally the metric product (splits along $\xi$ and $\ker\eta$).

\begin{lemma}
For a weak nearly Sasakian manifold $M^{\,2n+1}({f},Q,\xi,\eta,g)$ we obtain
\begin{eqnarray*}
 (\nabla_X\,h)\xi = -h(h-{f})X ,\qquad
 (\nabla_X{f})\,\xi =  -{f}(h-{f}) X .
\end{eqnarray*}
Moreover, if \eqref{E-nS-10} is true, then
\begin{eqnarray*}
 && h\,{f} + {f}\,h = -2\,\widetilde Q ,\quad
 h\,Q = Q\,h \quad (h\ {\rm commutes\ with}\ Q) , \\
 && h^2{f}={f} h^2,\quad h{f}^2={f}^2 h,\quad h^2{f}^2={f}^2 h^2.
\end{eqnarray*}
\end{lemma}

\begin{proposition}
Let a weak nearly Sasakian manifold
satisfy \eqref{E-nS-10}--\eqref{E-nS-04c}, then
\begin{equation*}
  g(R_{\,\xi, Z}\,{f} X,{f} Y) = 0,\quad {\rm hence,}\ \ker\eta\ {\rm is\ a\ curvature\ invariant\ distribution}.
\end{equation*}
\end{proposition}

\begin{lemma}\label{L-nS-04}
For a weak nearly Sasakian manifold $M^{\,2n+1}({f},Q,\xi,\eta,g)$ with conditions \eqref{E-nS-10} and \eqref{E-nS-04c}, we obtain
\begin{eqnarray}
\label{E-3.23d}
 & R_{\,\xi, X}Y = -(\nabla_X\,(h-{f})\,)\,Y, \\
\label{E-3.24d}
 & (\nabla_X\,(h-{f}))Y =  g( (h-{f})^2 X, Y)\,\xi - \eta(Y)\,(h-{f})^2 X , \\
\label{E-3.25d}
 & {\rm Ric}\,(\xi, Z) = -\eta(Z)\,({\rm tr}\,(h^2 + \widetilde Q) - 2\,n).
\end{eqnarray}
In particular, ${\rm tr}(h^2+ \widetilde Q) = const$, ${\rm Ric}\,(\xi, \xi) = const\ge0$ and
 $\nabla_\xi\,h = \nabla_\xi\,{f} = {f} h + \widetilde Q$.
By \eqref{E-3.23d}--\eqref{E-3.24d}, we get \
 $g(R_{\,\xi, X}Y, Z) =  \eta(Y)\,g( (h-{f})^2 X , Z) - \eta(Z)\,g( (h-{f})^2 X, Y)$.
\end{lemma}

\begin{proposition}
For a weak nearly Sasakian manifold with the property \eqref{E-nS-10}, the equality $h=0$ holds if and only if the manifold is Sasakian.
\end{proposition}

\begin{proposition}
For a weak nearly Sasakian manifold with conditions \eqref{E-nS-10} and \eqref{E-nS-04c},
the eigenvalues (and their multiplicities) of
the
self-adjoint operator $h^2$
are constant.
The spectrum of
$h^2$ has the~form~\eqref{E-nS-11b}:
  $Spec(h^2) = \{0, -\lambda_1^2,\ldots -\lambda_r^2\}$.
\end{proposition}

In particular, ${\rm tr}\,(h^2) = const\le0$, and by Lemma~\ref{L-nS-04}, ${\rm tr}\,Q = const>0$.

Denote by $[\xi]$ the 1-dimensional distribution generated by $\xi$,
and by $D_0$ a smooth distribution of the eigenvectors of $h^2$ with eigenvalue 0 orthogonal to $\xi$.
Denote by $D_i$ a smooth distribution of the eigenvectors of $h^2$ with eigenvalue $-\lambda^2_i$.
Note that the distributions $D_0$ and $D_i\ (i = 1,\ldots, r)$ belong to $\ker\eta$ and are ${f}$-invariant and $h$-invariant.
In particular, the eigenvalue $0$ has multiplicity $2p+1$ for some $p\ge0$.
If $X$ is a unit eigenvector of $h^2$ with eigenvalue $-\lambda^2_i$, then by \eqref{E-nS-01b} and \eqref{E-nS-01d},
$X, {f} X, hX$ and $h\,{f} X$ are nonzero eigenvectors of $h^2$ with eigenvalue $-\lambda^2_i$.
First, we show that weak nearly Sasakian manifolds satisfying \eqref{E-nS-10}--\eqref{E-nS-04c} have a foliated~structure.

\begin{theorem}
Let $M^{\,2n+1}({f},Q,\xi,\eta,g)$ be a weak nearly Sasakian
manifold with conditions \eqref{E-nS-10} and \eqref{E-nS-04c},
and let the spectrum of the self-adjoint operator $h^2$ have the form \eqref{E-nS-11b},
where the eigenvalue $0$ has multiplicity $2p+1$ for some integer $p\ge0$.
Then, the distribution $[\xi]\oplus D_0$ and
each distribution $[\xi]\oplus D_i\ (i = 1,\ldots, r)$ are integrable with totally geodesic leaves.
If $p>0$, then

\noindent\ \
$(a)$ the distribution $[\xi]\oplus D_1\oplus\ldots\oplus D_r$ is integrable and defines a (2n-2p+1)-dimensional
Riemannian foliation with totally geodesic~leaves;

\noindent\ \
$(b)$ the leaves of $[\xi]\oplus D_0$ are $(2p+1)$-dimensional Sasakian manifolds.
\end{theorem}

Next, we give some properties of the tensors ${f}$ and $h$,
and  two theorems characterizing Sasakian manifolds in the class of weak almost contact metric manifolds.

First, we consider weak almost contact metric manifolds with the condition \eqref{E-nS-10} and
characterize Sasakian manifolds in this class using the property \eqref{E-nS-Sas}.

\begin{theorem}
Let $M({f}, Q, \xi, \eta, g)$ be a weak almost contact metric manifold with conditions \eqref{E-nS-Sas} and \eqref{E-nS-10}.
Then $Q={\rm id}_{\,TM}$ and
$M({f}, \xi, \eta, g)$ is Sasakian.
\end{theorem}

Next, we generalize \cite[Proposition~3.1]{NDY-2018}.

\begin{proposition}
Let a weak nearly Sasakian manifold satisfy \eqref{E-nS-10}--\eqref{E-nS-04c}, then
\begin{eqnarray*}
 & (\nabla_X\,{f})Y = \eta(X)\,({f}  h Y +\widetilde Q Y) -\eta(Y)\,({f} h X + Q X) + g({f} hX + Q X, Y)\,\xi, \\
 & (\nabla_X\,h)Y = \eta(X)\,({f} h Y +\widetilde Q Y) -\eta(Y)\,h (h - {f}) X + g( h(h - {f}) X, Y)\,\xi,\\
 & (\nabla_X\,{f} h)Y = \eta(X)\,({f} h^2 Y - hY + \widetilde Q {f} Y )
  -\eta(Y)\,g({f} h^2 X - Q h X + 2\,\widetilde Q{f} X ) \\
 &\quad + g({f} h^2 X - hX + \widetilde Q hX, Y)\,\xi, \\
 & g((\nabla_{X}\,{f})Y, h Z) = -\eta(X)\,g( ({f} h^2 + \widetilde Q h) Z, Y) +\eta(Y)\,g( ({f} h^2 - h +\widetilde Q h) Z, X) .
\end{eqnarray*}
\end{proposition}

Recall,
that for any 2-form $\beta$ and 1-form $\eta$ we have
\begin{equation*}
 3\,(\eta\wedge\beta)(X,Y,Z) = \eta(X)\,\beta(Y,Z) + \eta(Y)\,\beta(Z,X) + \eta(Z)\,\beta(X,Y) .
\end{equation*}

\begin{proposition}
Let $\eta$ be a contact 1-form on a smooth manifold $M$ of dimension $2n+1> 5$
and $\Lambda^{p}(M)$ the vector bundle of differential $p$-forms on $M$.
Then, the operator
$\Upsilon_{d\eta} : \beta\in\Lambda^2(M) \to d\eta\wedge\beta\in\Lambda^4(M)$
is injective.
\end{proposition}

Using the above, we generalize Theorem~3.3 in \cite{NDY-2018}.

\begin{theorem}
Let a weak nearly Sasakian manifold $M({f}, Q, \xi, \eta, g)\ (\dim M>5)$
satisfy \eqref{E-nS-10}--\eqref{E-nS-04c}. Then $Q={\rm id}_{\,TM}$ and the structure $({f}, \xi, \eta, g)$ is Sasakian.
\end{theorem}


\section{Weak $\beta$-Kenmotsu manifolds}
\label{sec:08}
 \vglue-10pt
 \indent

Recall \cite{olszak1991normal} that the warped product $\mathbb{R}\times_{\sigma} \bar M$ (of $\mathbb{R}$ and a K\"{a}hler mani\-fold $(\bar M, \bar g)$)
with the metric $g=dt^2\oplus \sigma^2\,\bar g$ and the function $\sigma>0$ given on $(-\varepsilon, \varepsilon)$ and satisfying
\begin{equation}\label{E-f-warp}
 (\partial_t\,\sigma)/\sigma=\beta,
\end{equation}
admits a $\beta$-Kenmotsu structure $(\xi,\eta,{f})$, see
 \eqref{2.3-patra};
conversely, any point of a $\beta$-Kenmotsu manifold has a neighbourhood, which is a warped product $(-\varepsilon, \varepsilon)\times_{\sigma}\bar M$
of an interval and a K\"{a}hler mani\-fold, where $\sigma$ satisfies~\eqref{E-f-warp}.

\subsection{Geometry of weak $\beta$-Kenmotsu manifolds}
 \vglue-10pt
 \indent

Here, we show that a weak $\beta$-Kenmotsu manifold is locally the warped product $(-\varepsilon,\varepsilon)\times_\sigma\bar M$, where $(\partial_t\,\sigma)/\sigma=\beta\ne0$, and $\bar M$ has a weak Hermitian structure.

The following formulas are true for a weak $\beta$-Kenmotsu manifold:
\begin{equation*}
 \nabla_{X}\,\xi = \beta\{X -\eta(X)\,\xi\},\qquad
 d\eta(\xi,\,X)=0\qquad (X\in\mathfrak{X}_M).
\end{equation*}

\begin{proposition}\label{P-2.1}
A weak $\beta$-Kenmotsu manifold $M^{2n+1}({f},Q,\xi,\eta,g)$
with $n>1$ is a weak almost contact metric manifold satisfying
${N}^{\,(1)} {=} d\eta {=} 0$ and $3\,d\Phi = 2\beta\eta\wedge\Phi$.
\end{proposition}

The condition $d\beta\wedge\eta=0$ follows from \eqref{2.3-patra} if $\dim M>3$,
and it does not hold if $\dim M=3$. Indeed, by Proposition~\ref{P-2.1},
we get $0= 3\,d^2\Phi=2\,d\beta\wedge\eta\wedge\Phi$.

\begin{example}\label{Ex-2.1}\rm
Let $(\bar M,\bar g)$ be a Riemannian manifold.
A \textit{warped product} $\mathbb{R}\times_\sigma\bar M$
is the product $M=\mathbb{R}\times\bar M$ with the metric $g=dt^2\oplus \sigma^2(t)\,\bar g$, where $\sigma>0$ is a smooth function on $\mathbb{R}$.
Set $\xi=\partial_t$ and denote by ${\cal D}$ the distribution on $M$ orthogonal to $\xi$.
%
The~Levi-Civita connections, $\nabla$ of $g$ and $\bar\nabla$ of $\bar g$,
are related as follows:

\noindent\ \
(i) $\nabla_\xi\,\xi= 0$,
%
$\nabla_X\,\xi=\nabla_\xi\,X= \xi(\log\sigma)X$ for $X\in{\cal D}$.

\noindent\ \
(ii) $\pi_{1*}(\nabla_XY) = -g(X,Y)\,\xi(\log\sigma)$,
where $\pi_1: M \to\mathbb{R}$ is the orthoprojector.

\noindent\ \
(iii) $\pi_{2*}(\nabla_XY)$ is the lift of $\bar\nabla_XY$, where $\pi_2: M \to\bar M$ is the orthoprojector.

\smallskip\noindent
Submanifolds $\{t\}\times\bar M$ (called the \textit{fibers}) are totally umbilical, i.e., the Weingarten operator $A_\xi=-\nabla\xi$ on ${\cal D}$ is conformal with the factor $-\xi(\log\sigma)$, see (iii).
Note that $\sigma$ is constant along the fibers; thus, $X(\sigma)=\xi(\sigma)\,\eta(X)$ for all $X\in\mathfrak{X}_M$.

Let $(\bar M,\bar g,\bar f)$ be a weak K\"{a}hlerian manifold,
and
$\partial_t\sigma\ne0$.
Then the warped product $\mathbb{R}\times_\sigma\bar M$ is a weak $\beta$-Kenmotsu manifold with $\beta = (\partial_t\sigma)/\sigma$.
Indeed, define tensors,
\begin{equation}\label{E-struc}
 {f} = \left(\hskip-1mm
             \begin{array}{cc}
               0 & 0 \\
               0 & \bar f \\
             \end{array}
           \hskip-1mm\right),\
 Q = \left(\hskip-1mm
             \begin{array}{cc}
               0 & 0 \\
               0 & -\bar f^{\,2}
\\
             \end{array}
           \hskip-1mm\right),\
           \xi = \left(\hskip-1mm
            \begin{array}{c}
               \partial_t  \\
               0 \\
             \end{array}
           \hskip-2mm\right),\ \
        \eta = (dt, 0) ,\
           g = \left(\hskip-1mm
             \begin{array}{cc}
              dt^2 & 0 \\
               0 &
\sigma^2\bar g \\
             \end{array}
           \hskip-1mm\right)
\end{equation}
 on $M=\mathbb{R}\times\bar M$.
Note that
 $X(\beta)=0\ (X\perp\xi)$.
\newline
If $X,Y\in{\cal D}$, then
 $\pi_{2*}((\nabla_X{f})Y) =(\bar\nabla_X\bar{f})Y=0$
 and
 $\pi_{1*}((\nabla_X{f})Y) =-\beta\,g(X,{f} Y)$.
If~$X\in{\cal D}$ and $Y=\xi$, then
 $(\nabla_X{f})\xi=-{f}\nabla_X\,\xi=-\beta\,{f} X$.
If~$X=\xi$ and $Y\in{\cal D}$, then
$(\nabla_\xi{f})Y=\beta({f} Y)-{f}(\beta Y)=0$.
Also, we get $(\nabla_\xi{f})\xi=0$. By the above, \eqref{2.3-patra} is true.
%
\end{example}

\begin{theorem}\label{T-2.1}
Every point of a weak $\beta$-Kenmotsu manifold $M({f},Q, \xi,\eta,g)$ has a neighborhood $U$ that is a warped product $(-\varepsilon,\varepsilon)\times_\sigma\bar M$, where $(\partial_t\sigma)/\sigma=\beta$,
$(\bar M,\bar g,\bar f)$ is a weak K\"{a}hlerian manifold, and the structure $(\xi,\eta,{f},Q, g)$ is given on $U$ as \eqref{E-struc}.
Thus, the equality  $X(\beta)=0\ (X\perp\xi)$
is valid.
\end{theorem}

From Example~\ref{Ex-2.1} and Theorem~\ref{T-2.1}, we obtain the following generali\-zation of \cite[Theorem 2.3]{olszak1991normal}.
A weak almost contact metric manifold $M(f,Q,\xi,\eta,g)$
is a weak $\beta$-Kenmotsu manifold if and only if the following conditions are valid:
\newline\ \
-- the $\xi$-trajectories are geodesics and ${f}$ is $\xi$-invariant, i.e., $\pounds_\xi{f}=0$,
\newline\ \
-- the distribution $\ker\eta$ is integrable and defines a totally umbilical foliation ${\cal F}$ of codimension one,
whose leaves have constant mean curvature,
\newline\ \
-- a weak Hermitian structure $(\bar g, \bar f)$ induced on a leaf $\bar M\in{\cal F}$ is weak K\"{a}hlerian.

\begin{example}\rm
Let $(\bar M,\bar g,\bar f)$ be a weak K\"{a}hlerian manifold and
$\sigma(t)=c\,e^t\ (c=const\ne0)$ a function on a line $\mathbb{R}$. Then the warped product manifold $M=\mathbb{R}\times_\sigma\bar M$ has a weak almost contact metric structure which satisfies \eqref{2.3-patra} with $\beta\equiv1$.
\end{example}

\begin{lemma}
The following formulas hold on weak $\beta$-Kenmotsu manifolds:
\begin{eqnarray*}
 && R_{X, Y}\,\xi = (\xi(\beta)+\beta^2)(\eta(X)Y - \eta(Y)X)\quad (X,Y\in\mathfrak{X}_M),\\
 && \Ric^\sharp \xi = -2n(\xi(\beta)+\beta^2)\xi ,\\
 &&(\nabla_{\xi}\Ric^\sharp)X = -2\,\beta\Ric^\sharp X -2\,(\xi(\beta^2) +2\,n\,\beta^3)X +2\,(1-n)\,\xi(\beta^2)\,\eta(X)\,\xi,\\
 &&\xi(r) = -2\,\beta\,(r+2\,n(2\,n+1)\beta^2) -6\,n\,\xi(\beta^2),\quad X\in\mathfrak{X}_M;
\end{eqnarray*}
in particular, if $\beta=const$, then
$(\nabla_{\xi}\Ric^\sharp)X= -2\,\beta\Ric^\sharp X -4n\,\beta^3X$.
\end{lemma}

The following theorem and corollary generalize results in \cite{ghosh2019ricci} with $\beta\equiv1$.

\begin{theorem}
If a weak $\beta$-Kenmotsu manifold $M({f},Q,\xi,\eta,g)$ with $\beta=const\ne0$
satisfies $\nabla_{\xi}\Ric^\sharp=0$, then $(M,g)$ is an Einstein manifold with
$r=-2n(2n+1)\beta^2$.
\end{theorem}

\noindent\ \
A 3-dimensional Einstein manifold has constant sectional curvature. Thus, we~get

\begin{corollary}
Let $M^3({f},Q,\xi,\eta,g)$ be a
weak $\beta$-Kenmotsu mani\-fold with $\beta=const\ne0$.
If $\nabla_{\xi}\Ric^\sharp=0$, then $M$ has sectional curvature $-\beta^2$.
\end{corollary}

\subsection{$\eta$-Ricci solitons on weak $\beta$-Kenmotsu manifolds}
 \vglue-10pt
 \indent

Here, we study $\eta$-Ricci solitons on weak $\beta$-Kenmotsu manifolds.
We consider an $\eta$-Einstein weak $\beta$-Kenmotsu metric as an $\eta$-Ricci soliton
and characterize the Einstein metrics in such a wider class of metrics.

\begin{lemma}
Let $M^{2n+1}({f},Q,\xi,\eta,g)$ be a weak $\beta$-Kenmotsu manifold. If $g$ represents an $\eta$-Ricci soliton with the potential vector field $V$, then  for all $X\in\mathfrak{X}_M:$
\begin{eqnarray*}
 &(\pounds_V R)_{X,\,\xi}\,\xi = -2\,\xi(\beta)\Ric^\sharp X
 +\big( 4\,n\,\xi(\beta^3)
 +2\,\xi(\xi(\beta^2)) \big) X \\
 &-\big(8\,n\,\xi(\beta^3)
 +2(n+1)\,\xi(\xi(\beta^2)) \big)\,\eta(X)\,\xi;
\end{eqnarray*}
moreover, if
$\beta=const\ne0$, then $(\pounds_V R)_{X,\xi}\xi = 0$.
\end{lemma}

\begin{lemma}
On an $\eta$-Einstein weak $\beta$-Kenmotsu manifold $M^{2n+1}({f},Q,\xi,\eta,g)$, the expression of $\Ric^\sharp$ is the following:
\begin{equation*}
 {\rm Ric}^\sharp X = \big(\xi(\beta)+\beta^2+\frac{r}{2\,n}\big)\,X -\big((2n+1)(\xi(\beta)+\beta^2)+\frac{r}{2\,n}\big)\,\eta(X)\,\xi
 \quad (X\in\mathfrak{X}_M).
\end{equation*}
\end{lemma}

\begin{theorem}\label{thm3.1A}
Let $M^{2n+1}({f},Q,\xi,\eta,g)$, $n>1$, be an $\eta$-Einstein weak $\beta$-Kenmot\-su manifold with $\beta=const\ne0$.
If $g$ represents an $\eta$-Ricci soliton with the potential vector field $V$, then $(M,g)$ is an Einstein manifold with
$r=-2\,n(2\,n+1)\,\beta^2$.
\end{theorem}

\begin{definition}
\rm
 A vector field $X$ on a weak contact metric manifold
 is called
 a {\it weak contact vector field},
if
there exists a smooth function $\rho : M \rightarrow \mathbb{R}$ such~that
 $\pounds_{X}\,\eta=\rho\,\eta$,
and if $\rho=0$, then $X$ is said to be \textit{strict weak contact vector field}.
\end{definition}

We consider a weak $\beta$-Kenmotsu metric as an $\eta$-Ricci soliton, whose potential vector field $V$ is weak contact, or $V$ is collinear to $\xi$.
First, we derive the following.

\begin{lemma}\label{lem3.3}
 Let $M^{2n+1}({f},Q,\xi,\eta,g)$ be a weak $\beta$-Kenmotsu manifold.
If $g$ represents an $\eta$-Ricci soliton with a potential vector field $V$, then $\lambda+\mu=2n(\xi(\beta){+}\beta^2)$.
\end{lemma}

\begin{theorem}
Let $M^{2n+1}({f},Q,\xi,\eta,g)$, $n>1$, be a weak $\beta$-Kenmotsu manifold with
$\beta=const\ne0$ and $\xi(r)=0$.
If $g$ represents an $\eta$-Ricci soliton with a weak contact potential vector field $V$, then $V$ is strict
and $(M,g)$ is an Einstein manifold with scalar curvature $r=-2\,n(2n+1)\beta^2$.
\end{theorem}

\begin{theorem}
Let $M^{2n+1}({f},Q,\xi,\eta,g)$ be a weak $\beta$-Kenmotsu manifold with $\beta=const\ne0$ and $\xi(r)=0$. Suppose that $g$ represents an $\eta$-Ricci soliton with a non-zero
vector field $V$ collinear to $\xi$: $V=\delta\,\xi$ for a smooth function $\delta\ne0$ on~$M$.
Then $\delta=const$ and $(M,g)$ is an Einstein manifold with scalar curvature $-2n(2\,n + \delta(\beta-1))$.
\end{theorem}

In \cite{kenmotsu1972class}, K. Kenmotsu found the following necessary and sufficient condition for a Kenmotsu manifold $M$ to have constant
${f}$-holomorphic sectional curvature $H$:
\begin{eqnarray*}
 & 4R_{X,Y}Z=(H-3)\{g(Y,Z)X-g(X,Z)Y\} \\
 & +(H+1)\{\eta(X)\,\eta(Z)Y-\eta(Y)\,\eta(Z)X +\eta(Y)g(X,Z)\xi-\eta(X)g(Y,Z)\xi\\
 & +\,g(X,{f} Z){f} Y-g(Y,{f} Z){f} X+2g(X,{f} Y){f} Z\},\quad
 X,Y,Z\in\mathfrak{X}_M.
\end{eqnarray*}
The following corollary of Theorem~\ref{thm3.1A} illustrates Lemma~\ref{lem3.3}
and gives an example of a weak $\beta$-Kenmotsu manifold that admits an $\eta$-Ricci~soliton.

\begin{corollary}
Let a metric $g$ of a warped product $\mathbb{R}\times_{\sigma} N^{2n}$,
where $n>1$, $(\partial_t\sigma)/\sigma=\beta
=const\ne0$ and $N$ is a weak K\"{a}hlerian manifold, represent an $\eta$-Ricci soliton
with the potential vector field $V$.
Then $g$ has constant ${f}$-holomorphic sectional curvature $H=-\beta^2$.
\end{corollary}

We complete Theorem~\ref{thm3.1A} by studying the 3-dimensional case.

\begin{lemma}\label{L-3.3}
A
weak $\beta$-Kenmutsu manifold $M^3({f},Q,\xi,\eta,g)$ with $\beta=const\ne0$ represents an $\eta$-Einstein manifold.
\end{lemma}

Using Lemma~\ref{L-3.3} and
Theorem~\ref{thm3.1A}, we obtain the following.

\begin{corollary}
If a
weak $\beta$-Kenmutsu manifold
$M^{3}({f},Q,\xi,\eta,g)$ with 
$\beta=const\ne0$ and
$X(r) = 0\ (X\perp\xi)$
represents an $\eta$-Ricci soliton \eqref{1.1},
then $(M,g)$ has constant sectional curvature~$-\beta^2$.
\end{corollary}

\begin{lemma}
If a
weak $\beta$-Kenmutsu manifold
$M({f},Q,\xi,\eta,g)$ with $\beta=const\ne0$ represents an $\eta$-Ricci soliton \eqref{1.1}
with the potential vector field $V$,~then
\begin{eqnarray*}
&&(i)~2\,(\pounds_V\nabla)(X,\xi)=-\xi(r)\,(X-\eta(X)\xi),\quad {\rm where}\quad
\xi(r) = -2 \,\beta\,(\,r+6\,\beta^2),\\
&&(ii)~(r+6\,\beta^2)\,\big\{g(X,\pounds_{V}\xi) - \eta(X)\,\eta(\pounds_V\xi)\big\}
+ X(\xi(r)) + \xi(\xi(r))\,\eta(X) \notag\\
&&\qquad - 4\big\{\beta\,\xi(r)-2\,\beta^2\,(\lambda-2\,\beta^2)\big\}\,\eta(X) =0,\quad  X\in\mathfrak{X}_M.
\end{eqnarray*}
\end{lemma}


\begin{theorem}
If a weak $\beta$-Kenmutsu manifold $M^3({f},Q,\xi,\eta,g)$ with $\beta=const\ne0$ represents an $\eta$-Ricci soliton
with the potential vector field $V$, then $(M,g)$ has constant sectional curvature~$-\beta^2$.
\end{theorem}



\begin{thebibliography}{999}

\bibitem{blair2010riemannian}
{D.\ Blair}: \textit{Riemannian geometry of contact and symplectic manifolds}, Springer, 2010.

\bibitem{blair1974}
{D.\ Blair}: \textit{Almost contact manifolds with Killing structure tensors, II}.
J. Differential Geometry, {\bf 9} (1974), 577--582

\bibitem{blair1976}
{D.\ Blair, D.K.\ Showers, {and} Y.\ Komatu}: \textit{Nearly Sasakian manifolds}. {Kodai Math. Sem. Rep.}, {\bf 27},  (1976) 175--180.

\bibitem{bg2}\vskip-1.00mm
{Ch.P.\ Boyer, K.\ Galicki}: \textit{Sasakian Geometry}.
Oxford University Press, Oxford, 2008.

\bibitem{CZB-2022}
{A.M.\ Cherif, K.\ Zegga, and G.\ Beldjilali}: \textit{On the generalised Ricci solitons and Sasakian manifolds}.
\emph{Communications in Mathematics}, {\bf 30}, 1 (2022) 119--123.

\bibitem{cho2009ricci}
{J.\ Cho, and M.\ Kimura}: \textit{Ricci solitons and real hypersurfaces in a complex space form}, Tohoku Math. J. {\bf 61} No. 2, (2009) 205--212.

\bibitem{D-B-2021}
{S.\ Deshmukh, {and} O.\ Belova}: \textit{On Killing vector fields on Riemannian manifolds}, \emph{Mathe\-matics}, {\bf 9}, 259 (2021), 1--17

\bibitem{E-2005}
{H.\ Endo}: \textit{On the curvature tensor of nearly cosymplectic manifolds of constant $\varphi$-sectional curvature}.
Anal. \c{S}tiintifice Ale Univ. ``Al.I. Cuza" Ia\c{s}i Tomul LI, s.I, Matematic\u{a}, 2005, f.\,2

\bibitem{FP-2017}
{V.\ L.\,M. Fern\'{a}ndez {and} A. Prieto-Mart\'{i}n}: On $\eta$-Einstein para-$S$-manifolds.
Bull. Malays. Math. Sci. Soc. {\bf 40},  (2017) 1623--1637.

\bibitem{JKK-94}
{J.B.\ Jun, et al}:
\textit{On 3-dimensional almost contact metric manifolds}, Kyungpook Math. J. {1994}, {\bf 34}\,(2)  293--301.

\bibitem{ghosh2019ricci}
{A.\ Ghosh}: Ricci soliton and Ricci almost soliton within the framework of Kenmotsu manifold, Carpathian Math. Publ., {\bf 11}\,(1), 59--69, (2019).

\bibitem{G-D-2020}
{G.\ Ghosh {and} U.\,C\ De}: \textit{Generalized Ricci soliton on K-contact manifolds}, {Math. Sci. Appl. E-Notes}, {\bf 8} (2020), 165--169.

\bibitem{G-70}
{A.\ Gray}: \textit{Nearly K\"{a}hler manifolds}. J. Differ. Geom. {\bf 4}, 283--309 (1970).

\bibitem{H-2022}
{A.C.\ Herrera}: \textit{Parallel skew-symmetric tensors on 4-dimensional metric Lie algebras}.
Revista de la Uni\'{o}n Matem\'{a}tica Argentina, {\bf 65}, no. 2 (2023), 295--311.

\bibitem{kenmotsu1972class}
{K.\ Kenmotsu}: \textit{A class of almost contact Riemannian manifolds}, T\^{o}hoku Math. J., {\bf 24}, (1972) 93--103.

\bibitem{KKN}
{D.L.\ Kiran Kumar, {and} H.G.\ Nagaraja}:
\textit{Second order parallel tensor and Ricci solitons on genera\-lized $(k;\mu)$-space forms}.
Mathematical Advances in Pure and Applied Sciences, 2019, {\bf 2}, No. 1, 1--7.


\bibitem{NDY-2018}
{A.D.\ Nicola, G.\ Dileo, {and} I.\ Yudin}: \textit{On nearly Sasakian and nearly cosymplectic manifolds}.
Annali di Matematica Pura ed Applicata, {\bf 197}, (2018) 127--138.

\bibitem{olszak1991normal}
{Z.\ Olszak}: \textit{Normal locally conformal almost cosymplectic manifolds}, Publ. Math. Debrecen, {\bf 39}\,(3-4) (1991), 315--323.

\bibitem{Ol-1980}
{Z.\ Olszak}: \textit{Five-dimensional nearly Sasakian manifolds}. {Tensor New Ser.} {1980}, {\bf 34}, 273--276.

\bibitem{OV-2024}
{\sc L.\ Ornea, {\rm and} M.\ Verbitsky}: \textit{Principles of Locally Conformally K\"{a}hler Geo\-metry},
Progress in Mathematics, {\bf 354}, Birkh\"{a}user, Cham, 2024.

\bibitem{RWo-2}
{V.\ Rovenski {and} R.\ Wolak}: \textit{New metric structures on $g$-foliations}. Indagationes Mathema\-ticae, {\bf 33}, (2022) 518--532.


\bibitem{R2021}
{V.\ Rovenski {and} P.\ Walczak}: \textit{Extrinsic Geometry of Foliations}, Springer, Cham, 2021.


\bibitem{rov-119}
{V.\ Rovenski {and} D.S.\ Patra}: \textit{On the rigidity of the Sasakian structure and characterization of cosymplectic manifolds}.
Differential Geometry and its Applications, {\bf 90}, (2023) 102043.

\bibitem{rov-126}
{V.\ Rovenski {and} D.S.\ Patra}: \textit{Weak $\beta$-Kenmotsu manifolds and $\eta$-Ricci solitons}. pp. 53--72.
 In: Rovenski, V., Walczak, P., Wolak, R. (eds). Differential Geometric Structures and Applications,
2023. Springer Proceedings in Mathematics \& Statistics, {\bf 440}. Springer, Cham. 

\bibitem{rov-117}
{V.\ Rovenski}: \textit{Generalized Ricci solitons and Einstein metrics on weak K-contact manifolds}.
Communications in Analysis and Mechanics, 2023, {\bf 15}, No. 2: 177--188.

\bibitem{rov-122}
{V.\ Rovenski}: \textit{Weak nearly Sasakian and weak nearly cosymplectic manifolds}. Mathematics, 2023, {\bf 11}\,(20), 4377.

\bibitem{rov-128}
{V.\ Rovenski}: \textit{On the splitting of weak nearly cosymplectic manifolds},
Differential Geometry and its Applications, {\bf 94}, (2024) 102142.

\bibitem{rov-129}
{V.\ Rovenski}: \textit{Foliated structure of weak nearly Sasakian manifolds},
Annali di Matematica Pura ed Applicata (2024). https://doi.org/10.1007/s10231-024-01459-7

\bibitem{rov-130}
{V.\ Rovenski}: \textit{Characterization of Sasakian manifolds},
Asian-European Journal of Mathematics, {\bf 17}, No. 03, 2450030 (2024) (15 pages).

\bibitem{r-2024}
{\sc V.\ Rovenski}: \textit{Weak quasi contact metric manifolds and new characteristics of~K-contact and Sasakian manifolds}.
ArXiv [Math. DG]: 2410.02752.

\bibitem{YK-1985}
{K.\ Yano, {and} M.\ Kon}: \textit{Structures on Manifolds}, Vol. {\bf 3} of Series in Pure Math. World Scientific Publ. Co., Singapore, 1985.

\end{thebibliography}
\end{document}